\documentclass{elsart}

\usepackage{pstricks}
\usepackage{pstricks-add}
\usepackage{pst-node}
\usepackage{pst-key}
\usepackage{graphicx}
\usepackage{algorithmic}
\usepackage{algorithm}

\usepackage[utf8]{inputenc}

\usepackage{amsfonts}

\usepackage[centertags]{amsmath} 
\usepackage{amssymb}

\usepackage{amscd}
\usepackage{epsf}
\usepackage{epsfig}
\usepackage{cite}
\usepackage[section]{placeins}
\usepackage[english]{babel}

\def\:{\mkern 1.2mu \colon}

\newcommand{\simplex}[2]{\ensuremath{#1^{(#2)}}}
\newcommand{\simplexx}[1]{\ensuremath{#1}}
\newcommand{\dual}[1]{\ensuremath{#1^{*}}}

\newcommand{\mframe}[1]{\ensuremath{\texttt{fr}(#1)}}

\begin{document}

\begin{frontmatter}


\title{Ascending and descending regions of a discrete Morse function}


\author{Gregor Jerše}
\address{Institute of Mathematics, Physics and Mechanics\\ 
Ljubljana, Jadranska 19\\
Slovenia}
\ead{gregor.jerse@fmf.uni-lj.si}

\author{Neža Mramor Kosta}

\thanks{Both authors were partially funded by the Research Agency of the 
Republic of Slovenia, grant no. P1-0292}

\address{
University of Ljubljana, Faculty of Computer and Information Science\\
and Institute of Mathematics, Physics and Mechanics\\ 
Ljubljana, Jadranska 19\\
Slovenia}
\ead{neza.mramor@fri.uni-lj.si}

\begin{abstract}
We present an algorithm which produces a decomposition of a regular
cellular complex with a discrete Morse function analogous to the
Morse-Smale decomposition of a smooth manifold with respect to a
smooth Morse function. The advantage of our algorithm compared to
similar existing results is that it works, at least theoretically, in any
dimension. Practically, there are dimensional restrictions due to the
size of cellular complexes of higher dimensions, though. We prove that
the algorithm is correct in the sense that it always produces a
decomposition into descending and ascending regions of the critical
cells in a finite number of steps, and that, after a finite number of
subdivisions, all the regions are topological discs. The efficiency of
the algorithm is discussed and its performance on several examples is
demonstrated.

\end{abstract}

\end{frontmatter}

\newpage

%\maketitle

\section{Introduction and motivation}

Given a smooth manifold $M$ with a Morse function $F$ defined on it,
classical Morse theory \cite{milnor1963morse},
\cite{matsanxi2002ubka} is a powerful tool for investigating the 
behavior of the function $F$, as well as the topological properties of $M$. 
The critical points of $F$ together with their indices determine  
a handlebody decomposition of the domain $M$. A different decomposition of $M$, carrying 
topological as well as geometric information, is obtained from 
the stable and unstable manifolds of the critical points. The intersections 
of these are regions where the function behavior is uniform in the sense that 
all gradient paths have the same asymptotic behavior. If the function is Morse-Smale, that is, 
if the stable and unstable 
manifolds of the critical points intersect transversely, then these regions form the Morse-Smale CW-complex. 
Such a decomposition enables a 
thorough understanding of the domain $M$, as well as of the gradient flow of the function $F$.

\begin{figure}
\centerline{\includegraphics[width=8cm]{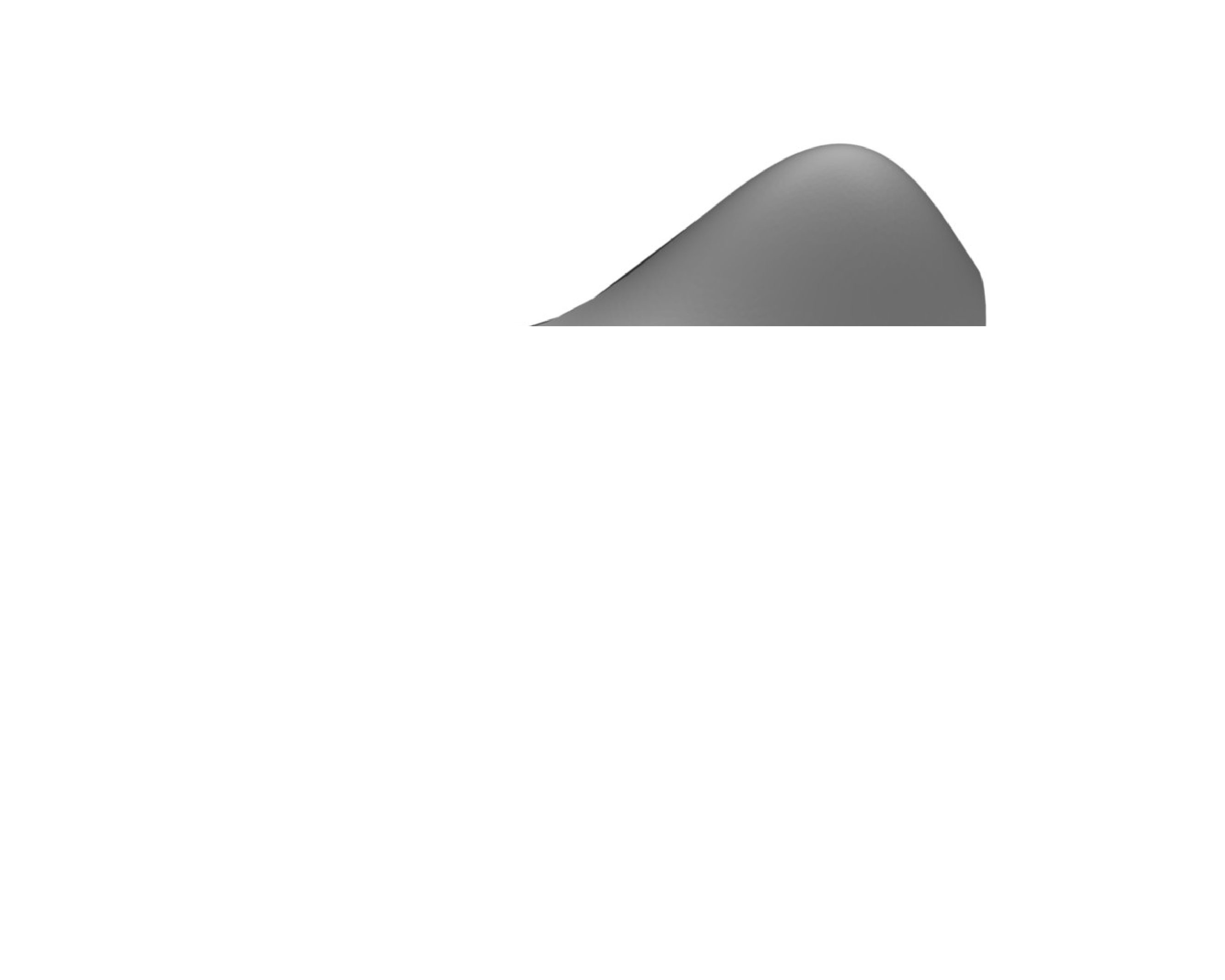}}
\caption {A model of a geographical terrain}
\label{fig:3D_terrain}
\end{figure}

Consider for example the function 
shown on figure \ref{fig:3D_terrain}. 
Assume that the function 
models a geographic terrain, and that our task is to find a path starting at a point $(x,y)$ in the  
descending disk of one of the maxima and ending in our preferred maximum $(x_0,y_0)$, which is traced 
optimally with respect to some criterion, for example height variation. 
This can be reconstructed from the descending and ascending disks, but these are  
computationally expensive. The problem is even more difficult if, instead of the function $f(x,y)$, 
only a sample of points on the surface is given. What 
we need, is a discrete approximation of the ascending and descending disks of $F$. 
Once this is given, the required path can be constructed so that it first reaches a 
saddle in the common boundary of both maxima, and then ascends towards the preferred one.

In this paper we propose an algorithm based on the discrete Morse theory of Forman 
\cite{formandiscrete}, \cite{forman01users} to solve this problem. Discrete Morse theory is a PL analogue 
of classical smooth 
Morse theory which has gained wide popularity and is used in topological data analysis
\cite{LLT04}, \cite{HKN}, \cite{lewiner-lopes-tavares-expmath03},  \cite{tbs} as well as in addressing purely 
theoretical topological and combinatorial problems, as for example in \cite{KC}, \cite{MY}, \cite{G}. 

In order to use discrete Morse theory for analyzing data, the initial data are first extended to a discrete 
Morse function on a triangulation (or, more generally, regular cellular decomposition) of the domain. 
This can be done using for example the algorithm of 
\cite{KKM05} (which is used in our implementations). 
The algorithm of this paper provides a further step by constructing the descending and ascending regions of the 
discrete Morse function and, from these, the discrete Morse-Smale complex.

The construction of discrete descending regions of a critical cell is 
motivated by the definition of descending disks in the 
smooth case: as unions of $V$-paths, which are the discrete analogue of
gradient paths, starting in the critical cell.
Discrete ascending
regions are constructed using the same procedure on the dual $V$-paths in the dual complex $K^*$. 
The resulting decomposition into descending regions is similar
to the smooth case in top dimension where the obtained regions are disjoint 
topological disks. In lower dimensions they are not necessarily disjoint, and not always disks, 
since, unlike gradient paths in the smooth case, $V$-paths from different critical cells or even from 
the same critical cell ending in the boundary of a specific critical cell $s$ can merge before reaching 
this boundary. 
We show, though, that after a finite number of subdivisions such merges can be 
eliminated, so that, as in the smooth case, the descending regions and the ascending regions form two families of 
disjoint topological disks. 
The regions obtained typically do not intersect transversely. This is not a major problem in the 
discrete case, though. The intersections nevertheless 
provide a discrete Morse-Smale decomposition of the cells of the domain into regions where, like in 
the continuous case, all $V$-paths come from and tend towards the same two critical cells. 

In \cite{G}, Theorem 3.1 the connection between discrete and 
smooth Morse theory is described in detail. It is 
shown that a smooth Morse function $F$ on a closed Riemannian manifold $M$ can be approximated by a 
discrete Morse function $f$ on a $C^1$ triangulation of $M$ so that critical points of $F$ correspond to 
critical cells of $f$. If, in addition, the function $F$ is Morse-Smale, 
the $V$-paths connecting any pair of critical cells $p$ and $q$ such that $\dim(p) = \dim(q) + 1$
are in bijection with integral curves of
the gradient vector field of $F$ between the corresponding critical points.

Several other algorithms for constructing discrete ascending and descending disks, as well as a discrete   
Morse-Smale complex exist. A classical example of an application of such a decomposition is the watershed 
segmentation algorithm in 
digital image analysis \cite{HPRCV,BM}. 
In \cite{EHZ03} a discrete Morse-Smale decomposition 
of a surface was constructed and applied to geographic surface modeling,
in \cite{cazals-chazal-lewiner-socg03}, using a different approach, such a decomposition was used 
in molecular modelling. In \cite{vams} 3D Morse-Smale complexes are used in volumetric data analysis. In \cite{edelsbrunner03morsesmale} an 
algorithm for constructing a quasi Morse-Smale complex, a combinatorial analogue of a 3D Morse-Smale 
complex, is given and in \cite{GNPBH} an extension of this algorithm is applied to several well known data sets. 
An overview of computational Morse theory with a large number of additional references and applications can be 
found in \cite{Z}. All existing algorithms are restricted to 2 or 3 independent variables, though, 
while real-world data often depends on more than 3 variables.

Our algorithm works, theoretically, in any dimension, and we have 
implementations in OCAML, C++ and C\# that accept simplicial or, more generally, polyhedral complexes
without dimensional restrictions as input.
In practice, however, the implementations can be used for simplicial
complexes of dimension 6 or less, since the number of cells in a
complex grows exponentially with the dimension, and the 
complexity of the algorithm depends on the number of cells. One of the applications
described in the last section uses a 4-dimensional data set.

An additional advantage of discrete Morse theory, in comparison to other approaches, is an efficient and elegant 
mechanism for dealing with noise which produces, after canceling pairs of critical cells, 
a simplified function in the sense that the number of critical elements is reduced. This is similar to 
handle sliding in computational Morse theory of \cite{edelsbrunner03morsesmale} and \cite{EHZ03}, and allows the use of 
persistence \cite{796607}, \cite{Zomorodian05computingpersistent}
%\cite{CZCG}
which has proven to be extremely useful in topological analysis of data \cite{DBLP:journals/tvcg/BremerEHP04}, 
\cite{CL}, \cite{GNPBH}.

The paper is organized as follows. In Section \ref{theory} we give a
short overview of discrete Morse theory and introduce notation. In
Section \ref{algoritem} we describe the construction of 
descending and ascending regions of critical cells, and discuss the  
complexity and the properties of the resulting decomposition. Finally, Section \ref{primeri} 
contains several demonstrations and applications of our algorithm. We first give a short description 
of 
an application of our algorithm to qualitative data analysis in artificial intelligence presented in \cite{QR07}. 
Next, an analysis of the response of a mechanical system consisting of a cart and a rod attached to its top 
is described. Our third example concerns the problem, presented in the beginning of this section. 
A procedure for constructing an optimal path (with respect to some optimization criterion) using the 
decomposition into descending regions is described, which is well defined also in the case, where the descending 
and ascending regions are not disks. Finally, 
we use our decomposition to construct a 
macroeconomic model on a 4-dimensional data set. This is part of an ongoing project with the goal to find and algorithm 
for controlling the key parameters of the model which contribute to long term economic growth of a country.

Part of this work was the result of a cooperation with the Laboratory
for artificial intelligence at the Faculty of Computer and Information
Science, University of Ljubljana. We would like to thank the members
of this lab for the permission to use data gathered in their experiments. Our thanks 
also go to the anonymous referees for valuable comments which have helped to improve this paper.

\section{Discrete Morse theory}\label{theory}

A discrete Morse function on a regular cellular complex $K$ associates
to every cell $\sigma\in K$ a real number $f(\sigma)$ such that in
most cases $f$ increases with increasing dimension except, possibly, in
one direction. More precisely, for every cell
$\simplex{\tau}{p} \in K$ of dimension $p$ the number of its codimension 1 faces
with values of $f$ greater than $f(\sigma)$ is at most one, and also
the number of its codimension 1 cofaces with values of $f$ smaller than $f(\sigma)$
is at most $1$. That is:

$$b(\tau) = \#\{\simplex{\nu}{p-1} \, | \, \simplexx{\nu}<\simplexx{\tau},  f(\simplexx{\nu}) \ge f(\simplexx{\tau})  \} \le 1$$
$$a(\tau) = \#\{\simplex{\sigma}{p+1} \, | \, \simplexx{\tau} < \simplexx{\sigma}, f(\simplexx{\tau}) \ge f(\simplexx{\sigma})  \} \le 1.$$

\begin{figure}[ht]
  \centering
	\psset{unit=2.5cm}
	\begin{pspicture}(0,0)(2,1)

        \rput[bl](0,0){\rnode{A}{\simplexx{\nu}}}
        \rput[bc](0.75,0){\rnode{B}{\simplexx{\tau'}}}
	\rput[tl](0.3,0.6){\rnode{C}{\simplexx{\tau}}}
        \rput[cc](0.85,0.6){\rnode{D}{\simplexx{\sigma}}}
	\psline(0.1,0.1)(1.6, 0.1)
        \psline(0.1,0.1)(0.85, 1.1)
        \psline(1.6,0.1)(0.85, 1.1)
        %arrows
        \psline[linewidth=0.03]{->}(0.1, 0.1)(0.4 , 0.5)			
        \psline[linewidth=0.03]{->}(0.43 , 0.5)(0.75, 0.5)
        \psline[linewidth=0.03]{->}(0.85, 0.5)(0.85 , 0.15)
        \psline[linewidth=0.03]{->}(0.85 , 0.095)(0.2 , 0.095)			
	\end{pspicture} 
  \caption{A discrete Morse function can not have both $a(\tau)$ and $b(\tau)$ equal to $1$}
  \label{figure:dm_basic}
\end{figure}
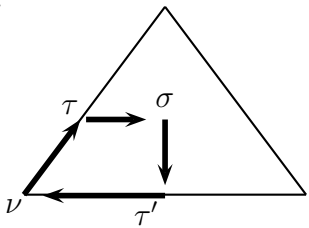

The numbers $a(\tau)$ and
$b(\tau)$ can not both be one. Since, if there exists a face $\nu< \tau$ 
such that $f(\nu)\geq f(\tau)$ as well as a coface $\sigma> \tau$ such that 
$f(\sigma)<f(\tau)$, then for any $\tau'<\sigma$ such that $\nu<\tau'$, 
as in figure \ref{figure:dm_basic} (where the arrows indicate the direction 
of function descent) it follows that 
$$ f(\tau') > f(\nu) \ge f(\tau)\ge f(\sigma) > f(\tau') $$ 
which is a contradiction.  

Because of this, the cells of $K$ split into three subsets $K=A\cup B\cup C$ (in the terminology of \cite{KKM05}). 
If  both  $a(\simplexx{\tau})=b(\simplexx{\tau})=0$ then $\tau$ is a {\em critical cell} 
of index $p=\dim \tau$, and $C$ is the set of all critical cells. 
If not, either $a(\tau)=1$ or $b(\tau)=1$, and $\tau$ is a
{\em regular cell}. A regular cell \simplexx{\tau} belongs to 
the set $A$ if $a(\simplexx{\tau})=1$ and to $B$ if 
$b(\simplexx{\tau})=1$.

For each $\tau\in A$ there exists precisely one $\sigma\in B$
such that $\tau < \sigma$ and $f(\tau)\geq f(\sigma)$. The map $V:\ A
\to B$ which associates $\simplex{\tau}{p}$ to
$\simplex{\sigma}{p+1}$ is called a {\em discrete gradient vector field}, and 
points in the direction of steepest descent from the cell $\tau$.

It is convenient to imagine a discrete gradient vector field not as a
map but as a collection of pairs $(\simplex{\tau}{p},
\simplex{\sigma}{p+1})$ such that 
$V(\simplexx{\tau})=\simplexx{\sigma}$, or arrows from
$\simplexx{\tau}$ to $\simplexx{\sigma}$, in the direction of function
descent. This representation of a discrete gradient vector field is
used in all figures in this paper. The cells of $K$ which are unpaired
are precisely the critical cells.

In classical Morse theory the gradient vector field points in the
direction of steepest ascent of the function $f$, and its gradient paths (traced in the direction 
opposite to the gradient flow) are the paths of steepest descent. 
The discrete analogue is a
$V$-path, which is a sequence of cells 
$$ \simplex{\tau_0}{p} <
\simplex{\sigma_0}{p+1} > \simplex{\tau_1}{p} < \simplex{\sigma_1}{p+1}
... \simplex{\tau_r}{p} < \simplex{\sigma_r}{p+1} >
\simplex{\tau_{r+1}}{p} $$ 
such that $V(\simplexx{\tau_{i}}) = \simplexx{\sigma_i}$ and 
$\simplexx{\tau_{i+1}}$ is a face of $\simplexx{\sigma_i}$ different 
from 
 $\simplexx{\tau_i}$ for $i\in \{1,..,r\}$.

As in the classical case, function values decrease along a
$V$-path. This implies that a $V$-path of a discrete gradient vector field $V$ arising from a discrete Morse function 
can not be closed, that is, the discrete gradient vector field can not contain any cycles. 
Forman showed that this is characteristic for 
discrete gradient vector fields: a pairing $V$ arises from a
discrete gradient vector field of a discrete Morse function on $K$ if and only
if it contains no nontrivial closed $V$-paths.

In classical Morse theory, two critical points $p_1$ and $p_2$ of a
Morse function $F$ on $M$ of indices differing by $1$ can be
\textit{canceled} when there exists exactly one gradient path between
them. That is, using handle sliding, $F$ can be smoothly deformed to a function $F'$ which
has exactly the same critical points as $F$ except $p_1$ and $p_2$
that become regular. This is true also in the discrete case. If the
critical cells $\simplexx{\sigma}$ and $\simplexx{\tau}$ are connected
by exactly one $V$-path starting in any cell of
$\partial\simplexx{\sigma}$ and ending in $\simplexx{\tau}$, i.e. $$
\simplex{\sigma}{p+1} = \simplex{\sigma_{-1}}{p+1} >
\simplex{\tau_0}{p} < \simplex{\sigma_0}{p+1} > \simplex{\tau_1}{p} <
\simplex{\sigma_1}{p+1} \ldots  \simplex{\tau_r}{p} <
\simplex{\sigma_r}{p+1} > \simplex{\tau_{r+1}}{p} = \simplex{\tau}{p}
$$ then the pair $\simplexx{\sigma},\simplexx{\tau}$ can be
cancelled. This is done simply by modifying the discrete gradient
vector field $V$ along the connecting paths to 
$$V'(\simplexx{\tau_i}) = \simplexx{\sigma_{i-1}} $$ for
$i\in{0,...,r+1}$. This procedure can also be described as switching
the direction of the arrows along the connecting $V$-path. Since there
are no other $V$-paths connecting
$\simplex{\tau_0}{p}<\partial\simplexx{\sigma}$ and
$\simplex{\tau}{p}$ this can not create a non-trivial closed $V'$-path
and so, according to Forman's characterization, the modified $V'$ is the
discrete gradient vector field of some Morse function for which the
cells \simplexx{\tau} and \simplexx{\sigma} are no longer critical.

\section{Algorithm}
\label{algoritem}

In this section we present our algorithm which, on the basis of a
discrete Morse function $f$, produces two decompositions of a regular
cellular complex $K$ with $|K|$ a manifold. First, a decomposition 
of $K$ into {\em discrete descending regions} of the critical cells of $f$ 
is constructed. Each of the
regions contains exactly one critical cell \simplex{s}{p}, the maximum
of the region, and a collection of regular cells of dimension less
than or equal to $p$.

In the second step we decompose, using the same procedure, the dual complex
\dual{K} into discrete descending regions of the function $-f$. This gives 
a decomposition of $K$ into {\em discrete ascending regions} of the critical 
cells.

\subsection{Constructing the descending regions}

In the smooth case, the unstable manifold of a critical point $c$ of $F$ is
the union of all gradient paths of the flow $-$grad $F$ starting close to the point $c$. In the
discrete case the direction of function descent is indicated by the
discrete gradient vector field $V$, which plays a similar role as the
gradient vector field in the smooth case. The discrete descending
region $D(\simplex{s}{p})$ of a critical cell $\simplex{s}{p}$ thus
contains all regular cells which appear in a $V$-path starting in
the boundary of $\simplex{s}{p}$. The descending region is constructed
in two steps. In the first step the {\em frame} which contains all cells
of maximal dimension is constructed, and in the second step the cells of 
lower dimension are added. 

The critical cells are first ordered by 
ascending dimension so that, when processing a critical cell of dimension $d$, the 
descending regions of all critical cells of dimension less than $d$ have already been constructed.
This will be very important in the second step of the algorithm.

Let \simplex{s}{p} be a critical cell. The frame \mframe{\simplexx{s}} consists of all regular 
$p$ and $(p-1)$ dimensional cells appearing in a $V$-path starting in the boundary of $\simplexx{s}$. 
This part of the algorithm is a simple breadth-first search through the regular $(p-1)$-faces of the included 
$p$-cells which terminates when no more regular $(p-1)$-cells belonging to the 
set $A$ (that is, no more $(p-1)$-dimensional arrow tails) can be found. 
Since a cell in the frame can typically be reached by several paths in $V$, a set structure is used 
for storing the cells to eliminate duplicates. 

\begin{figure}[ht]
  \centerline{\mbox{
		\input {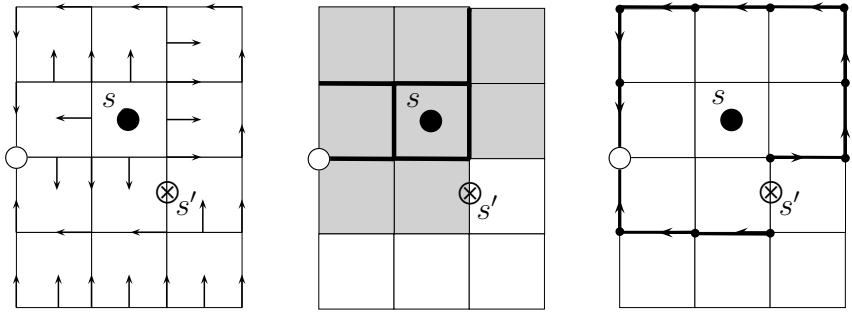}
  }}
  \caption{Building the frame of descending regions in a cellular complex}
  \label{figure_frame}
\end{figure}

This step of the algorithm is illustrated in figure \ref{figure_frame}. On the left, the discrete gradient 
vector field is shown, the descending frame of the maximum is in the middle 
and the descending frame of the saddle is on the right.

%\begin{algorithm}[ht!]
%\caption{Frame for the critical cell \simplex{s}{p}} 
%\label{framealg}
%\begin{algorithmic}
% \STATE let $q$ be an empty queue
% \STATE let $fr$ be an empty set
% \STATE add simplex \simplex{s}{p} to the queue $q$
% \STATE add simplex \simplex{s}{p} to the set $fr$
% \WHILE {queue q is not empty}
%   \STATE remove the first element from queue $q$
%   \STATE let $faces$ be the set of its regular $(p-1)$-dimensional faces
%   \FORALL {$f$ in $faces$}
%     \IF{$f$ is paired with a $p$-dimensional cell $g$}
%       \STATE add simplices $f$ and $g$ to the set $fr$
%       \STATE add simplex $g$ to the queue $q$
%     \ENDIF
%   \ENDFOR
% \ENDWHILE
% \RETURN $fr$
%\end{algorithmic}
%\end{algorithm}

The frame \mframe{\simplex{s}{p}} contains only $p$ and $p-1$ dimensional cells. 
The next step is to decide which additional faces of the cells in \mframe{\simplexx{s}} 
should be included.

Consider a $V$-pair $(\simplexx{\alpha}, \simplexx{\beta})$ incident to 
some simplex from \mframe{\simplex{s}{p}}. Since $V(\simplexx{\alpha}) = \simplexx{\beta}$, 
function values decrease from \simplexx{\alpha} towards \simplexx{\beta}. 
This implies that \simplexx{\beta} belongs to the descending region $D(\simplexx{s})$ 
when \simplexx{\alpha} does, so the pair $(\simplexx{\alpha}, \simplexx{\beta})$ is 
included into or excluded from the descending region together. 
It is included  
if all codimension one cofaces of $\alpha$ except $\beta$ have already been included. 
This is implemented 
by the following recursively used routine. Let $S$ be the set of all $(p+1)$-dimensional cofaces of 
\simplexx{\alpha}, different from $\beta$, and incident to some $p$-dimensional cell from \mframe{\simplexx{s}}. 
The routine checks if all elements of $S$ 
belong to $D(\simplexx{s})$. 
It ends when either a cell $\simplex{\gamma}{d+1} \in S$ is found that belongs to the descending region 
of some other critical cell of dimension less than $p$, in this case the pair (\simplexx{\alpha},\simplexx{\beta}) is not 
included, or all cells of $S$ have been included in $D(\simplexx{s})$, and in this case the 
pair $(\alpha,\beta)$ is also included. 

\begin{figure}[ht]
  \centerline{\mbox{
		\input {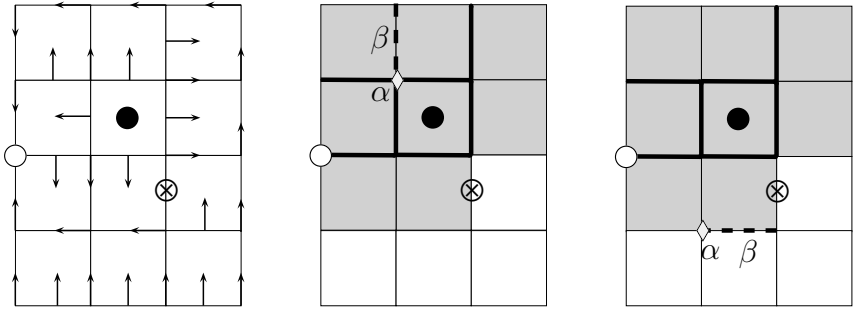}
  }}
  \caption{Completing the descending region}
  \label{figure_completing_region}
\end{figure}

On figure \ref{figure_completing_region} the left picture again shows the discrete gradient vector field. 
In the middle picture, part of the descending region of the maximum is already constructed (the 2-cells forming the frame 
are shaded, and the 1-cells that are already included are bold), and the pair $(\alpha,\beta)$ that is processed is the 
diamond-shaped point and the dashed edge. Since all edges that have this point as a 
face are already included in the descending disk, the pair $(\alpha,\beta)$ is added. 
On the right, the situation is different. The descending disk of the saddle has already been 
constructed and contains 
the dashed edge $\beta$. Since $\beta$ already belongs 
to the descending region of the saddle, the pair $(\alpha,\beta)$ will not be included.

\begin{prop}
The construction of the descending region of a critical cell ends in a 
finite number of steps.
\end{prop}
\begin{pf}
During frame construction we simply follow the direction of the
discrete gradient vector field $V$. Since $V$ does not contain
nontrivial closed paths and the complex $K$ is finite the frame is 
completed in a finite number of steps.
 
The second step of the construction, where lower dimensional
cells are added, is recursive. The algorithm runs through all
$V$-pairs in the boundary of frame cells.
For each pair $(\simplex{\alpha}{d}, \simplex{\beta}{d+1})$ which
has not been processed it checks if all cofaces of $\alpha$ of dimension $(d+1)$ 
are already included in the descending region. If it encounters a coface which has not been
processed it checks that first. This implies that the algorithm
always proceeds in the direction opposite to $V$, so no
cycles can be created.

Since the number of pairs is finite, the algorithm ends in a finite
number of steps. Notice that the algorithm can start processing pairs
in any order, independent of dimension.
\end{pf}

\subsection{The boundary}

The algorithm described in the previous section constructs discrete
descending regions of a discrete Morse function on a complex $K$ with
$|K|$ a manifold without boundary. For manifolds with boundary, the
boundary must be processed separately.

Denote the discrete gradient vector field of the restriction of $f$ to the boundary by $V_{\partial K}$. 
A critical cell $\nu$ of $V_{\partial K}$ 
is 
\begin{enumerate}
	\item either a critical cell in $V$
	\item or is paired to a higher dimensional cell that does not belong to the $\partial K$; 
	such cells will be called {\em  boundary
critical cells}.
\end{enumerate}

In the first case no extra processing is needed, as $\nu$ is also a
critical cell in $V$ and its descending region has been
constructed in the previous step. 

In the second case, the discrete descending region of a boundary critical cell
\simplex{\nu}{p} is constructed in two stages. First, the
discrete descending region of the critical cell \simplex{\nu}{p} in
$\partial K$ is found using the algorithm of the
previous section. In the second stage, all cells of dimension $p$
in the discrete descending region of \simplex{\nu}{p}  in $\partial K$, that are paired by $V$ to a higher
dimensional cell $\beta$ in the interior, are considered. The cell $\beta$ 
is processed as a critical
cell and its descending region is constructed. The union of all regions
obtained in this way constitutes the descending region of the boundary
critical cell \simplexx{\nu}.

For example, consider $f(x,y) = x + y$ on a triangulation of the square
$[0,1]\times[0,1]$. The discrete gradient vector field
contains all
cells except a minimum on the boundary which is the only critical cell

\begin{figure}[ht!]
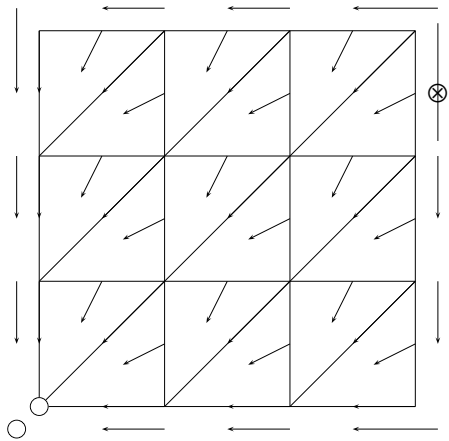

  \centerline{\mbox{
  		\psset{unit=5mm}
\begin{pspicture}(-0.2, -0.2)(10.2, 10.2)
  		\input {examples/boundary/graph.tex}
  		\input {examples/boundary/field.tex}
   		\input {examples/boundary/critical.tex}

  		\input {examples/boundary/square/field.tex}
  		\input {examples/boundary/square/critical.tex}
\end{pspicture}
  }}
  \caption{Discrete gradient vector field on the square and its boundary}
  \label{boundary_example_boundary}
\end{figure}

The discrete gradient vector fields of $f$ and $f|_{\partial K}$ are
shown on Figure \ref{boundary_example_boundary}. The maximum $\mu$ of
$f|_{\partial K}$ is a boundary critical $1$-cell. Its descending region in $\partial K$ 
is all of $\partial K$ except for the minimum. A boundary cell $\alpha$ which is paired 
with an interior cell $\beta$ in $V$ is included into the descending region of $\mu$ together 
with $\beta$ and its descending region. Finally, the descending region of $\mu$ is 
the entire square without the minimum.

%\newpage
\subsection{Ascending regions}

The ascending regions of the critical cells of $f$ on $K$ are obtained
from the dual complex \dual{K} with the dual the discrete gradient
vector field \dual{V}. 

Viewing $K$ as an abstract regular cell complex determined by its cells with specified dimension 
and the face relationship between them, 
the dual complex $K^*$ is obtained simply by reversing the dimensions and the face relationship on the cells. 
That is, a cell $\alpha^{(p)}\in K$ corresponds to a cell $\alpha^{*(n-p)}\in K^*$, where $n$ is the dimension 
of $K$, and if $\alpha$ is a face of $\beta$ in $K$, then $\alpha^*$ contains $\beta^*$ as its face in $K^*$. 
In the geometric realization $|K|$ of $K$, the vertices of $K^*$ correspond to barycenters of the $n$-cells of $K$, 
while a cell of higher dimension is spanned by a collection of vertices of $K^*$ such that the corresponding $n$-cells 
of $K$ have a nonempty intersection. If $|K|$ is a PL-manifold without boundary, then $|K^*|=|K|$. If 
$|K|$ is a PL-manifold with boundary, then $K^*$ is not a cell decomposition of $|K|$ any more, since the 
cells of $K^*$ which are dual to boundary cells of $K$ do not have a sufficient number of faces (for example, 
a $1$-cell of this type has only $1$ face). For our purpose $K^*$ is just a combinatorial object, though, so this 
geometric property does not represent a problem.

Also the dual discrete vector field $V^*$ is obtained simply by reversing the arrows in $V$. 
That is, two cells \dual{\alpha} and
\dual{\beta} are paired in \dual{V} if the cells
\simplexx{\beta} and \simplexx{\alpha} are paired in the original discrete
gradient vector field $V$. Thus, a cell \simplexx{\alpha} in $K$
is critical of index $p$ if and only of the cell \dual{\simplexx{\alpha}} is
critical of index $(n-p)$ in \dual{K}. 

\begin{figure}[ht]
  \centerline{\mbox{
		% Generated with LaTeXDraw 2.0.0
% Mon Sep 15 18:08:06 CEST 2008
% \usepackage[usenames,dvipsnames]{pstricks}
% \usepackage{epsfig}
% \usepackage{pst-grad} % For gradients
% \usepackage{pst-plot} % For axes
\scalebox{1} % Change this value to rescale the drawing.
{
\begin{pspicture}(0,-1.88)(6.005,1.865)
\definecolor{color1818b}{rgb}{0.8588235294117647,0.8588235294117647,0.8588235294117647}
\pspolygon[linewidth=0.0035277777,fillstyle=solid,fillcolor=color1818b](0.0017638779,0.4982361)(0.0017638889,1.4982361)(1.0017639,1.4982361)(1.0017639,0.49823612)
\pspolygon[linewidth=0.0035277777,fillstyle=solid](1.0017639,0.49823612)(1.0017639,1.4982361)(2.0017638,1.4982361)(2.0017638,0.49823612)
\pspolygon[linewidth=0.0035277777,fillstyle=solid,fillcolor=color1818b](0.0017638779,-0.50176394)(0.0017638889,0.49823612)(1.0017639,0.49823612)(1.0017639,-0.5017639)
\pspolygon[linewidth=0.0035277777,fillstyle=solid](1.0017639,-0.5017639)(1.0017639,0.49823612)(2.0017638,0.49823612)(2.0017638,-0.5017639)
\pspolygon[linewidth=0.0035277777,fillstyle=solid,fillcolor=color1818b](0.0017638779,-1.5017639)(0.0017638889,-0.5017639)(1.0017639,-0.5017639)(1.0017639,-1.5017639)
\pspolygon[linewidth=0.0035277777,fillstyle=solid](1.0017639,-1.5017639)(1.0017639,-0.5017639)(2.0017638,-0.5017639)(2.0017638,-1.5017639)
\psdots[dotsize=0.3,dotstyle=otimes](0.5317638,-0.5076389)
\psline[linewidth=0.02cm,arrowsize=0.05291667cm 2.0,arrowlength=1.4,arrowinset=0.4]{->}(1.0,1.02)(0.56,1.02)
\psline[linewidth=0.02cm,arrowsize=0.05291667cm 2.0,arrowlength=1.4,arrowinset=0.4]{->}(0.52,0.5)(0.52,0.04)
\psline[linewidth=0.02cm,arrowsize=0.05291667cm 2.0,arrowlength=1.4,arrowinset=0.4]{->}(0.52,-1.5)(0.52,-1.06)
\psdots[dotsize=0.3,dotstyle=oplus](1.491764,1.012361)
\pspolygon[linewidth=0.0035277777,fillstyle=solid](3.581764,0.4982361)(3.581764,1.4982361)(4.5817637,1.4982361)(4.5817637,0.49823612)
\pspolygon[linewidth=0.0035277777,fillstyle=solid](4.5817637,0.4982361)(4.5817637,1.4982361)(5.5817637,1.4982361)(5.5817637,0.49823612)
\pspolygon[linewidth=0.0035277777,fillstyle=solid](3.581764,-0.50176394)(3.581764,0.49823612)(4.5817637,0.49823612)(4.5817637,-0.5017639)
\pspolygon[linewidth=0.0035277777,fillstyle=solid](4.5817637,-0.50176394)(4.5817637,0.49823612)(5.5817637,0.49823612)(5.5817637,-0.5017639)
\pspolygon[linewidth=0.0035277777,fillstyle=solid](3.581764,-1.5017639)(3.581764,-0.5017639)(4.5817637,-0.5017639)(4.5817637,-1.5017639)
\pspolygon[linewidth=0.0035277777,fillstyle=solid](4.5817637,-1.5017639)(4.5817637,-0.5017639)(5.5817637,-0.5017639)(5.5817637,-1.5017639)
\psline[linewidth=0.04cm,arrowsize=0.05291667cm 2.0,arrowlength=1.4,arrowinset=0.4]{->}(4.14,1.0)(4.6,1.0)
\psline[linewidth=0.04cm,arrowsize=0.05291667cm 2.0,arrowlength=1.4,arrowinset=0.4]{->}(4.12,-1.02)(4.12,-1.48)
\psline[linewidth=0.04cm,arrowsize=0.05291667cm 2.0,arrowlength=1.4,arrowinset=0.4]{->}(4.12,0.0)(4.12,0.44)
\pspolygon[linewidth=0.01,linestyle=dashed,dash=0.16cm 0.16cm](4.1217637,-0.001763916)(4.1217637,0.9982361)(5.1217637,0.9982361)(5.1217637,-0.0017638889)
\pspolygon[linewidth=0.01,linestyle=dashed,dash=0.16cm 0.16cm](4.1217637,-1.0017639)(4.1217637,-0.0017638889)(5.1217637,-0.0017638889)(5.1217637,-1.0017639)
\psline[linewidth=0.01cm,linestyle=dashed,dash=0.16cm 0.16cm](4.12,0.0)(3.08,0.0)
\psline[linewidth=0.01cm,linestyle=dashed,dash=0.16cm 0.16cm](6.0,1.0)(5.12,1.0)
\psline[linewidth=0.01cm,linestyle=dashed,dash=0.16cm 0.16cm](4.12,-1.0)(3.08,-1.0)
\psline[linewidth=0.01cm,linestyle=dashed,dash=0.16cm 0.16cm](6.0,0.0)(5.12,0.0)
\psline[linewidth=0.04cm](4.12,-1.86)(4.12,-1.0)
\psline[linewidth=0.01cm,linestyle=dashed,dash=0.16cm 0.16cm](5.12,1.0)(5.12,1.86)
\psline[linewidth=0.01cm,linestyle=dashed,dash=0.16cm 0.16cm](6.0,-1.0)(5.12,-1.0)
\psline[linewidth=0.01cm,linestyle=dashed,dash=0.16cm 0.16cm](5.12,-1.86)(5.12,-1.0)
\psline[linewidth=0.01cm,linestyle=dashed,dash=0.16cm 0.16cm](4.12,1.0)(3.08,1.0)
\psline[linewidth=0.01cm,linestyle=dashed,dash=0.16cm 0.16cm](4.12,1.0)(4.12,1.86)
\psline[linewidth=0.04cm](4.12,1.02)(4.12,0.02)
\psline[linewidth=0.04cm](4.12,1.0)(5.12,1.0)
\psdots[dotsize=0.12,fillstyle=solid,fillcolor=black,dotstyle=o](4.111764,-1.0076392)
\psdots[dotsize=0.12,fillstyle=solid,fillcolor=black,dotstyle=o](4.111764,-0.02763916)
\psdots[dotsize=0.12,fillstyle=solid,fillcolor=black,dotstyle=o](4.111764,0.99236083)
\psline[linewidth=0.04cm](0.0,0.5)(1.0,0.5)
\psline[linewidth=0.04cm](0.0,-1.5)(1.0,-1.5)
\psline[linewidth=0.04cm](1.0,1.5)(1.0,0.5)
\psdots[dotsize=0.3,fillstyle=solid,dotstyle=o](5.111764,1.0123612)
\psline[linewidth=0.02cm](4.12,-1.06)(4.12,0.2)
\psdots[dotsize=0.3,fillstyle=solid,dotstyle=o](4.131764,-0.5076389)
\psdots[dotsize=0.3,dotstyle=otimes](4.131764,-0.5076389)
\end{pspicture} 
}
  }}
  \caption{Ascending region of a saddle}
  \label{figure_dual_disc}
\end{figure}
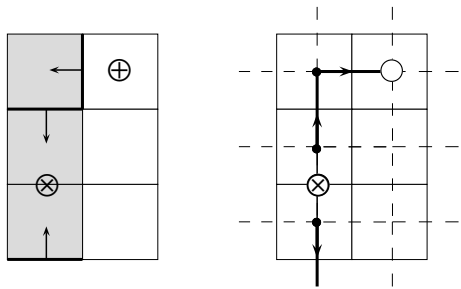

The discrete ascending region of a
critical cell \simplexx{\alpha} of $K$ is the dual of
the discrete descending region of the critical cell \dual{\simplexx{\alpha}} 
of \dual{K} (see figure \ref{figure_dual_disc}).

%\begin{algorithm}[ht!]
%\caption{Constructing the ascending area for the critical cell \simplex{s}{p}} 
%\label{ascendingalg}
%\begin{algorithmic}
%\STATE let $\dual{K}$ be the dual complex of complex $K$
% \STATE let $\dual{V}$  \texttt{be the empty set of pairs of cells}
% \FORALL {pairs $\{\simplexx{\alpha}, \simplexx{\beta}\}$ from $V$}
%  \STATE add pair $\{\dual{\simplexx{\beta}}, \dual{\simplexx{\alpha}}\}$ to $\dual{V}$
% \ENDFOR
%
% \STATE let $desc$ be the descending area for the cell $\dual{\simplexx{s}}$ in $(\dual{K}, \dual{V})$
% \STATE let $asc$ be the empty set of cells
% \FORALL{cells $\dual{p}$ in $desc$}
%  \STATE add cell $p$ to the set asc
% \ENDFOR 
% \RETURN $asc$
%\end{algorithmic}
%\end{algorithm}

%\FloatBarrier

\subsection{Consistency}

The descending and ascending disks of a smooth Morse function $f$ on a manifold $M$ 
form a decomposition of $M$ into topological disks. If the function is Morse-Smale, their intersections 
form the Morse-Smale complex on $M$. 
In the discrete case, the descending and ascending regions obtained are always disks only in the top 
dimension, while in lower dimensions 
they might not be disks. In addition, they might not be disjoint, that is, a cell can belong 
to the descending or ascending disk of more than one critical cell. We first prove that the descending 
regions (and therefore also the ascending regions) of all critical cells form a covering of the 
underlying PL manifold.

\begin{prop}
Let $K$ be a regular $n$-dimensional cellular complex such that $|K|$
is a manifold without boundary and let $V$ be the discrete gradient
vector field of a discrete Morse function on $K$. Then every
regular cell $\simplexx{\alpha}
\in K$ is contained in the discrete descending region of a critical cell
$\simplexx{\beta}$.

This is true also for manifolds with boundary, if the discrete
gradient vector field of $f$ on $K$ is an extension of a discrete
gradient vector field of $f\mid_{\partial K}$.
\end{prop}

\begin{pf}
We first prove that every regular top dimensional cell \simplex{\alpha}{n} is included in 
the frame of some top dimensional critical cell. A regular cell of top 
dimension is paired either to a cell \simplex{\gamma_0}{n-1} which is a 
face of at least one
$n$-dimensional cell \simplexx{\alpha_1} other than
\simplexx{\alpha_0}, or it is paired to a boundary critical cell \simplex{\nu}{n-1}
of top dimension in $\partial K$. In the second case, \simplex{\alpha}{n}
belongs to the descending disk of \simplex{\nu}{n-1}.
In the first case, if \simplexx{\alpha_1} is critical, then
\simplexx{\alpha_0} lies in its descending region, and if it is regular, we repeat
the same process for 
\simplexx{\alpha_1}. Since the discrete gradient vector field $V$ does
not have cycles and $K$ is a finite cellular complex this process must
end in a finite number of steps at a critical cell
$\simplex{\alpha_i}{n}$ or a boundary critical cell, and \simplexx{\alpha_0} belongs to its
descending disk.

Assume now that descending disks of all critical cells have 
been built and that there exists a regular pair $(\simplex{\alpha}{p}, \simplex{\beta}{p+1})$ that 
is not included in any of them. 
 
Since every cell is a face of some top-dimensional cell, \simplexx{\alpha} was 
processed during the construction of the descending region of some cell $\simplexx{\sigma}$. 
Since it was not included, there exists a coface $\simplex{\gamma_{1}}{p+1}$ of 
$\simplexx{\alpha}$, such that $\simplexx{\gamma_1}$ is included in the descending region 
$D(\simplex{\sigma_1}{n_1})$, where $n_1 < n$. 
In the construction of the descending region $D(\simplex{\sigma_1}{n_1})$ also $\simplexx{\alpha}$ was 
processed since it is also a face of $\simplexx{\gamma_0}$, and not included. So  
there must exist a coface $\simplex{\gamma_{2}}{p+1}$ of $\simplexx{\alpha}$ which belongs to 
some other descending region $D(\simplex{\sigma_2}{n_2})$, where $n_2 < n_1$. 
This argument can be repeated indefinitely, and we obtain an infinite sequence of critical cells 
$\simplexx{\sigma_0}, \simplexx{\sigma_1}, \simplexx{\sigma_2} \ldots $ with strictly decreasing positive 
dimensions which is a contradiction. So every regular pair $(\simplex{\alpha}{p}, \simplex{\beta}{p+1})$ is included in the descending region of some critical cell.

\end{pf}

In smooth Morse theory the stable and unstable manifolds of a critical
point are constructed from gradient paths, that is, integral curves of 
the gradient vector field, which begin
or end at this point, respectively. An important property of integral curves of 
a smooth vector field is
that they are pairwise disjoint - two different integral curves can be
arbitrarily close but do not merge. The discrete analogue of a gradient path is a $V$-path, and 
since there is no
such thing as \lq arbitrarily close\rq\/ in the discrete world, a $V$-path 
typically splits into multiple paths, 
and several $V$-paths can merge into one path. Due to this, 
discrete descending and ascending regions of 
different critical cells might not be disjoint, since $V$-paths beginning in 
different critical cells might merge. They also might not be 
topological disks, except in the highest dimension.

\begin{prop}
\label{descending_regions_of_maximal_critical_cells_are_discs}
The descending region of a critical cell of maximal dimension (i.e. a
maximum) collapses to the maximal cell. The same is true for ascending regions 
of a critical cell of dimension $0$ (i.e. a minimum).
\end{prop}

\begin{pf}
For every cell $\sigma$ of maximal dimension, its pair is a cell of
codimension $1$ which is either a boundary critical cell or it is the
face of precisely one other cell of maximal dimension. Because of this
there exists precisely one $V$-path starting at some (possibly
boundary) critical cell $s$ of maximal index (i.e. a maximum) leading
to $\sigma$. Since a discrete vector field has no cycles, every such 
$V$-path leading to a cell in the boundary of the descending region of $s$ 
determines a sequence of elementary collapses starting in the boundary of the 
descending region and ending in $s$. A finite number of such $V$-paths 
determines the necessary elementary collapses required to collapse 
the whole descending region to the cell $s$. 
\end{pf}

\begin{figure}[ht]
  \centerline{\mbox{
    \includegraphics[width=4.00cm,height=4.00cm]{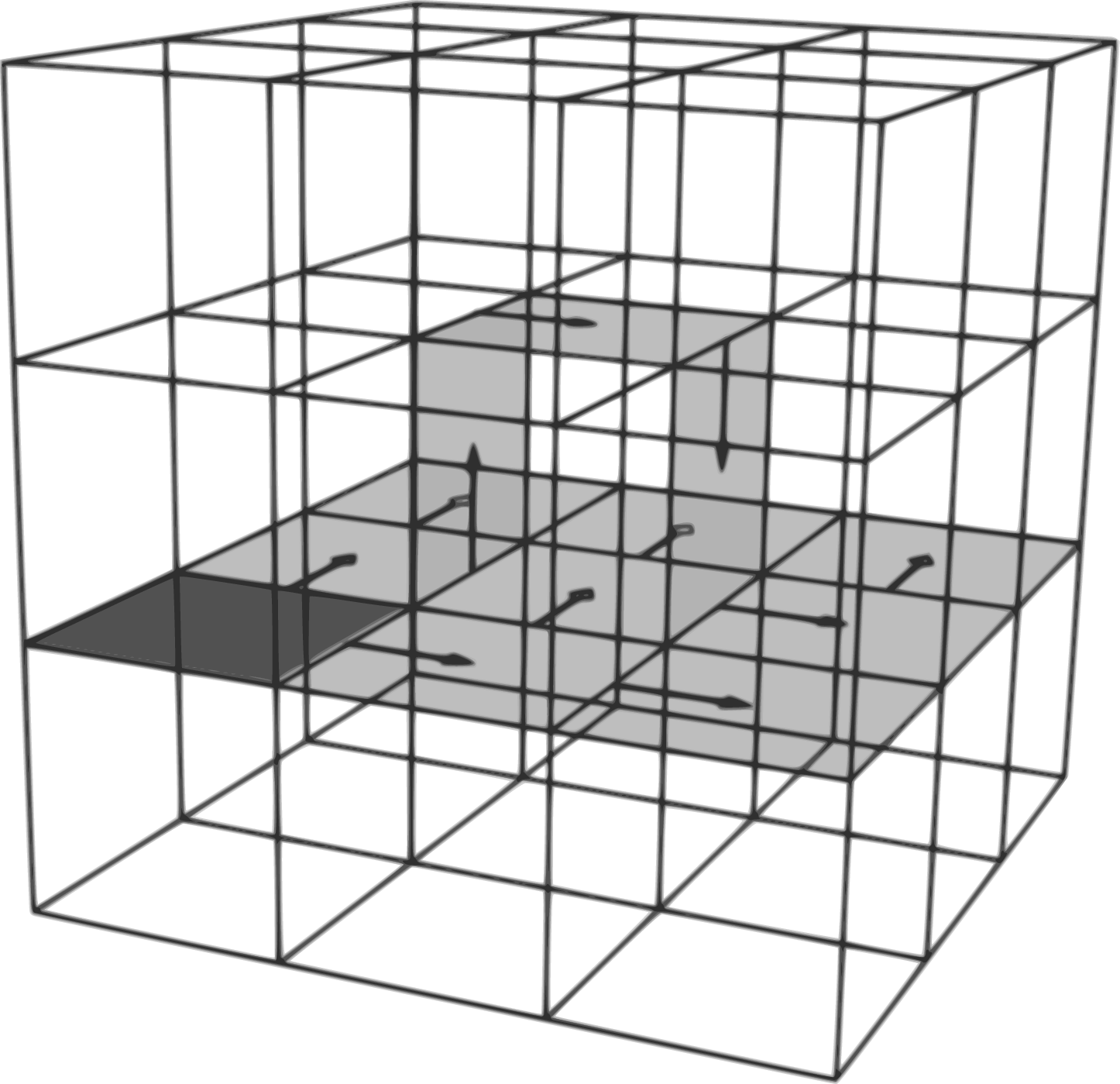}
  }}
  \caption{Discrete descending region of a $2$-saddle with splitting and merging $V$-paths}
  \label{figure_cube}
\end{figure}

For descending regions of critical cells of index smaller than
maximal, i.e. saddles, this is not necessarily true, since $V$-paths
not containing cells of maximal dimension can merge. Figure
\ref{figure_cube} shows a section of a 3-dimensional cubical
complex where $V$-paths in the descending region of a saddle split and
merge again, so that the union is not a disk.

It is possible, however, to modify the decomposition of the cellular
complex and the discrete gradient vector field so that the merging
point of descending paths is pushed further towards the boundary of
the descending region, while the discrete vector field outside the
open star of the merging point is left unaffected. The idea is similar
to splitting multiple saddles in \cite{EHZ03}.

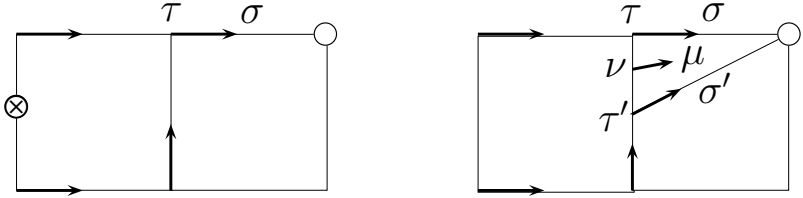
\begin{figure}[ht]
  \centerline{\mbox{
 		\scalebox{1.3}{% Generated with LaTeXDraw 2.0.0
% Mon Sep 15 18:46:52 CEST 2008
% \usepackage[usenames,dvipsnames]{pstricks}
% \usepackage{epsfig}
% \usepackage{pst-grad} % For gradients
% \usepackage{pst-plot} % For axes
\scalebox{1} % Change this value to rescale the drawing.
{
\begin{pspicture}(0,-0.97703123)(8.16,1.0170312)
\psframe[linewidth=0.0040,dimen=outer,fillstyle=solid](1.72,0.6429688)(0.12,-0.95703125)
\psframe[linewidth=0.0040,dimen=outer,fillstyle=solid](3.3,0.6429688)(1.7,-0.95703125)
\psline[linewidth=0.03cm,arrowsize=0.05291667cm 2.0,arrowlength=1.4,arrowinset=0.4]{->}(0.12,0.6429688)(0.8,0.6429688)
\psline[linewidth=0.03cm,arrowsize=0.05291667cm 2.0,arrowlength=1.4,arrowinset=0.4]{->}(0.12,-0.95703125)(0.8,-0.95703125)
\psline[linewidth=0.03cm,arrowsize=0.05291667cm 2.0,arrowlength=1.4,arrowinset=0.4]{->}(1.7,0.6429688)(2.38,0.6429688)
\psline[linewidth=0.03cm,arrowsize=0.05291667cm 2.0,arrowlength=1.4,arrowinset=0.4]{->}(1.7,-0.95703125)(1.7,-0.27703124)
\psframe[linewidth=0.0040,dimen=outer,fillstyle=solid](6.44,0.62296873)(4.84,-0.97703123)
\psframe[linewidth=0.0040,dimen=outer,fillstyle=solid](8.02,0.6429688)(6.42,-0.95703125)
\psline[linewidth=0.03cm,arrowsize=0.05291667cm 2.0,arrowlength=1.4,arrowinset=0.4]{->}(4.84,0.6429688)(5.52,0.6429688)
\psline[linewidth=0.03cm,arrowsize=0.05291667cm 2.0,arrowlength=1.4,arrowinset=0.4]{->}(4.84,-0.95703125)(5.52,-0.95703125)
\psline[linewidth=0.03cm,arrowsize=0.05291667cm 2.0,arrowlength=1.4,arrowinset=0.4]{->}(6.42,0.6429688)(7.1,0.6429688)
\psline[linewidth=0.03cm,arrowsize=0.05291667cm 2.0,arrowlength=1.4,arrowinset=0.4]{->}(6.42,-0.95703125)(6.42,-0.47703126)
\psline[linewidth=0.0040cm](6.4,-0.17703125)(8.02,0.6429688)
\psline[linewidth=0.03cm,arrowsize=0.05291667cm 2.0,arrowlength=1.4,arrowinset=0.4]{->}(6.42,-0.17703125)(6.92,0.08296875)
\psline[linewidth=0.03cm,arrowsize=0.05291667cm 2.0,arrowlength=1.4,arrowinset=0.4]{->}(6.42,0.28296876)(6.84,0.36296874)
\usefont{T1}{ptm}{m}{n}
\rput(1.7085937,0.87296873){$\tau$}
\usefont{T1}{ptm}{m}{n}
\rput(2.535,0.85296875){$\sigma$}
\usefont{T1}{ptm}{m}{n}
\rput(6.4085937,0.85296875){$\tau$}
\usefont{T1}{ptm}{m}{n}
\rput(6.233125,-0.18703125){$\tau'$}
\usefont{T1}{ptm}{m}{n}
\rput(7.255,0.87296873){$\sigma$}
\usefont{T1}{ptm}{m}{n}
\rput(7.2660937,0.09296875){$\sigma'$}
\usefont{T1}{ptm}{m}{n}
\rput(6.248125,0.29296875){$\nu$}
\usefont{T1}{ptm}{m}{n}
\rput(7.0376563,0.39296874){$\mu$}
\psdots[dotsize=0.24,fillstyle=solid,dotstyle=o](3.28,0.6429688)
\psdots[dotsize=0.24,fillstyle=solid,dotstyle=o](8.02,0.6429688)
\psdots[dotsize=0.24,fillstyle=solid,dotstyle=o](0.12,-0.09703125)
\psdots[dotsize=0.24,dotstyle=otimes](0.12,-0.09703125)
\end{pspicture} 
}}
  }}
	\caption{One step of pushing a merging point towards the boundary}
	\label{figure_notcontractible2d}
\end{figure}

The left example on figure
\ref{figure_notcontractible2d} shows a section of a 2-dimensional
simplicial complex where the descending region of a critical cell of
index $1$ is topologically a 1-sphere. The right side of this figure
shows the modification which pushes the merging point further along
descending values of $f$.  The figure shows only one step of the pushing procedure, which can be applied as many times as necessary.

Here is a description of the required modification of a regular
cellular complex in the general case. Assume that two $V$-paths
$$\ldots \simplex{\tau_{i-1}}{p-1}<\simplex{\sigma_{i-1}}{p}>
\simplex{\tau_i}{p-1}=\simplexx{\tau}<
\simplex{\sigma}{p}>\simplex{\tau_{i+1}}{p-1},\ldots$$ 
$$\ldots
\tau'^{(p-1)}_{j-1}<\sigma'^{(p)}_{j-1}>
\tau'^{(p-1)}_{j}=\simplexx{\tau}< \sigma^{(p)}> \tau'^{(p-1)}_{j+1},
\ldots$$ 
merge in the cell $\simplex{\tau}{p-1}$ which is paired with
$\simplex{\sigma}{p}$. Let $S$ be $\texttt{star}(\tau)$ and let $p_1$
%= \simplex{\sigma_{i-1}}{p} = \simplex{\alpha_1}{p}, \simplex{\beta_1}{p+1}, \simplex{\alpha_2}{p}, \simplex{\beta_2}{p+1}, \ldots \simplex{\beta_k}{p+1},\sigma$ 
and $p_2$ % =\simplex{\sigma_{i-1}}{p} = \alpha'^{(p)}_1,
%\beta'^{(p+1)}_1, \alpha'^{(p)}_2, \beta'^{(p+1)}_2, \ldots
%\beta'^{(p+1)}_k, \sigma'_{j-1}$
be two disjoint paths in $S$ from the
cell $\simplexx{\sigma_{i-1}}$ to the cells $\sigma$ and
$\sigma'_{j-1}$ respectively.

The required modification consists of replacing $\simplex{\tau}{p-1}$
by two parallel copies $\simplex{\tau}{p-1}$ and
$\simplex{\tau'}{p-1}$, and $\simplex{\sigma}{p}$ by two copies
$\simplex{\sigma}{p}$ and $\simplex{\sigma'}{p}$ with
$\simplexx{\partial \sigma'}$ containing the same cells as
$\simplexx{\partial\sigma}$, except for $\tau$ which is replaced by
$\tau'$. In addition, two new cells are added: a cell
$\simplex{\nu}{p}$ which has cells \simplex{\tau}{p-1} and
\simplex{\tau'}{p-1} in its boundary, and a cell \simplex{\mu}{p+1}
with $\nu, \sigma$ and $\sigma'$ in its boundary. We modify the star $S$
in the following way
\begin{itemize}
 \item in the boundary of all $p$-cells in the path $p_1$
 (except $\sigma$) the cell $\tau$ is replaced by $\tau'$,
\item in the boundary of the last $(p+1)$-cell in the 
path $p_1$ $\sigma$ is replaced by $\sigma'$,
\item $\nu$ is added to the boundary of the first $(p+1)$-cell in the path 
$p_2$,
\item $\tau'$ and $\nu$ are added to the boundary of all other cells in $S$.

\end{itemize}

The modified field $V'$ coincides with $V$ on both copies of
$\simplexx{\tau}$, $\simplexx{\sigma}$ and their boundaries, and it
pairs $\simplexx{\nu}$ with $\simplexx{\mu}$. The
new field $V'$ now has two $V'$-paths, which possibly merge in some cell
$\tau_{i+k}=\tau'_{j+l}$ (or are disjoint). So we have
pushed the merging point further along the discrete gradient
vector field $V$. Also note that this step does not affect $V$-paths originating from critical 
cells in 
dimension $p-1$ or lower and does not create any additional merge points for critical paths 
originating from critical cells of dimension $p$.

\begin{prop}
Let $K$ be a finite cellular complex without boundary and let $V$ be 
a discrete gradient vector field on $K$. Then there exists a finite sequence of steps described 
above which modify $K$ into $K'$ and $V$ into $V'$ in such a way that the set of critical cells remains 
unchanged, and that the descending regions of all critical cells with respect to $V'$ are topological 
disks. 
\end{prop}

\begin{pf}
The proof goes  by ascending dimension of the cells of the complex. Since there is nothing to do for 
$0$-dimensional critical cells, we start with critical cells of dimension $1$. 
Let $\simplex{\tau_1}{0}$ 
be the cell where $V$-paths originating from a critical $1$-cell \simplexx{s} merge and 
let $D$ be the descending region of $s$. 
Since $K$ is finite, a finite sequence of steps described above 
pushes the cell $\simplex{\tau_1}{0}$ into the boundary of $D$. We repeat this process 
for all such merges in the descending region of $s$.
Let $K''$ be the resulting complex, $V''$ the resulting gradient vector field and  
$D''$ the discrete descending region of $s$ in $K'$. 
Since $D''$ does not contain its boundary, and since all merging points have been pushed into 
the boundary we have pushed the cells where $V''$-paths originating 
from \simplexx{s} merge away from $D''$. Now the same argument as in the proposition 
\ref{descending_regions_of_maximal_critical_cells_are_discs} can be used to show that $D''$  
collapses to the critical cell $s$ and is therefore a topological disk. Also, all $V$-paths 
coming from critical cells different from $s$ have been pushed off $D''$, so $D''$ does 
not intersect any other descending regions.

Since all modifications take place in the interior of the descending disk of $s$, all other 
descending regions remain unchanged. And since no additional merge points have been 
created in the process it follows that, after repeating this procedure for all critical $1$-cells, 
all descending regions originating from critical cells of dimension 
one are topological disks.

Precisely the same arguments work for critical cells of higher dimension. Assuming that 
descending regions for critical cells of dimension less than $p$ 
have already been processed and are therefore topological disks, the descending 
regions for each critical cells of dimension $p$ is processed in the same way as above. Since
all merging points are pushed to the boundary and since no other descending regions are 
affected, all descending regions of critical cells of dimension $p$ are 
transformed one by one into disks. 
Moreover, since the pushing step does not affect the descending regions 
of critical cells of dimension less than $p$, descending regions of all critical cells of 
dimension up to $p$ are now topological disks.

\end{pf}

\subsection{Complexity}

Let $m$ be the total number of cells of $K$, $n$ the dimension of the
complex $K$, $m_d$ the number of cells of dimension $d$, $r_d$ the average number of
codimension $1$ faces (if $K$ is a simplicial complex $r_d=d+1)$, $p_d$ the average number of 
codimension $1$ cofaces and $p_{\max}=\max(p_0, p_1, \ldots , p_n)$. Let $c_t$ be
the total number of critical cells of the discrete Morse function $f$
on $K$, and $c_d$ the number of critical cells of dimension $d$.

The frame of a $d$-dimensional critical cell \simplex{s}{d} consists
of $d$ and $(d-1)$-dimensional cells that belong to $V$-paths starting
in the boundary of \simplexx{s}. For each $d$-dimensional cell in the
frame, the set of its $(d-1)$-dimensional faces has to be found, and
for each face its pair in $V$ has to be found.
Since faces and their pairs can be found in constant time,
this takes at most $\alpha \times r_d \times m_d$ operations for a suitable
constant $\alpha$. 

Building all frames of critical cells thus takes $\sum_{i=1}^{n} \alpha
\times m_i \times r_i $ operations. In the case of simplicial
complexes, the number of operations is at most $(n+1) \times \alpha \times m \times
n$. The complexity of this step is therefore
$O(n^2\times m)$.

To complete the descending regions requires determining whether a
regular cell \simplexx{q} that is incident to the region of some
descending disk is included in this region. First, its V-pair
\simplexx{p} has to be found which can be done in
constant time. Next, a list of all co-faces of the
lower dimensional cell of the pair is required which can also be found in constant
time. Each cell in this list is then processed in the same manner. 
On average it takes $p_j$ operations for each cell of dimension $j$, so we need at most 
$\sum_{i=1}^{d} \alpha''
\times p_i \times m_i $ operations to complete the descending region of one critical cell of dimension $d$, 
where $\alpha''$ is a constant. Since we can not use previous results when completing frames for other critical 
cells some cells are typically checked more than once. So completing a frame for all critical cells takes at most 
$c_t \sum_{i=1}^{d} \alpha''
\times p_{\max} \times m  = c_t \times p_{\max}\times m\times n$ operations.

For simplicial complexes the complexity of the algorithm is altogether
$O(nm\times (n+c_{t}p_{\max}))$. Note that the algorithm complexity increases exponentially with the
dimension of the complex, since the number of cells grows exponentially. 
In the special case where $K$ is the 
Delaunay triangulation on a set of points in $\mathbb{R}^n$, an upper bound for the total number of cells 
is $O(v^{\lceil n/2 \rceil})$ where $v$ is the number of points. In this case the algorithm has 
complexity $O(v^{\lceil n/2 \rceil}n \times (n+c_{t}p_{\max}))$.

\section{Examples}\label{primeri}

The algorithm presented in this paper has been applied to data sets from different domains. 
In this section we present some of them.

\subsection{The $QING$ algorithm}
\label{qing}
We first mention the algorithm $QING$ which is an application of our algorithm to AI from \cite{QR07}. 
The algorithm uses discrete Morse theory to reconstruct the critical cells of a discrete Morse
function obtained from the values of a sampled function. The descending disks obtained using an early implementation of 
this algorithm were used for the construction of a {\em qualitative graph}, that is, an undirected graph 
$G = (V_G, E_G)$ that represents 
connections between the critical cells. The set of vertices $V_G$ is the set of all critical cells. 
Two points $\simplexx{\alpha}$ and \simplexx{\beta} from $V_G$ with dim $\simplexx{\alpha}\geq$ dim $\simplexx{\beta}$ 
are connected if, $\simplexx{\beta}$ is in the boundary of the descending region of $\simplexx{\alpha}$.
This graph was then used in learning in a qualitative model. 

The $QING$ algorithm, together with the parametric discrete Morse theory of \cite{HKN} was used in a 
learning scheme presented in \cite{QR08} designed to teach a robot the concept of occlusion.

\subsection{Mechanical system response}
Our next example is obtained from a model constructed for the purpose 
of qualitative modeling of data in artificial intelligence. The model is taken from 
\cite{Jure_diploma}. We would like to thank the authors for permission to use
the data. The data represents measurements obtained from a system
consisting of a cart and a rod attached to its top in such a way that
it can fall either backwards or forwards. The cart is set on ice (or
other low-friction surface) and the rod is positioned upwards. As the 
rod begins to fall, the cart responds, accelerating in the opposite direction of the fall.

Our discrete Morse function models the response of the cart in terms of 
the acceleration $\ddot{x}$, depending on the rod angle $\varphi$ and the angular velocity $\dot{\varphi}$. 
The measured values of $\ddot{x}$ are given on a square grid. A discrete gradient vector
field was first constructed on a triangulation
of the grid using the algorithm of \cite{KKM05} producing, after cancellation, three
maxima, four saddles and two minima.

\begin{figure}[ht!]
  \centerline{\mbox{
    \includegraphics[height=10.00cm]{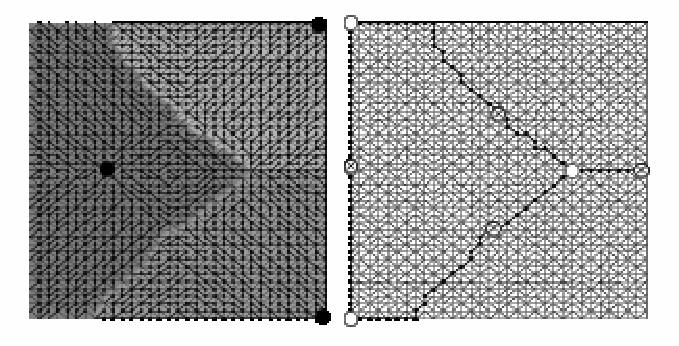}
  }}
 \caption{Discrete descending regions of maxima and saddles of the cart function}
 \label{cart_discs}
\end{figure}

On figure \ref{cart_discs} the discrete descending regions
of all three maxima are on the left, and discrete descending
regions of the four saddles are on the right. The
descending regions form a decomposition of the triangulated area into
disjoint discs. Note that the critical cell on the right side of the
image is a boundary saddle. Its descending area is correctly generated and separates the 
descending regions of the two maxima.

%\begin{figure}[ht!]
% \psset{unit=5mm}
% \centerline{\mbox{
% \begin{pspicture}(0, 0)(20, 10)
% 		\input {examples/cart/cart_maximum_discs.tex}
% 		\rput(11,0){
% 		  \input {examples/cart/cart_saddle_discs.tex}
% 		}  
% \end{pspicture}
% }}
% \caption{Discrete descending regions of maxima and saddles of the cart function}
% \label{cart_discs}
%\end{figure}

\begin{figure}[ht!]
 \centerline{\includegraphics[width=8cm]{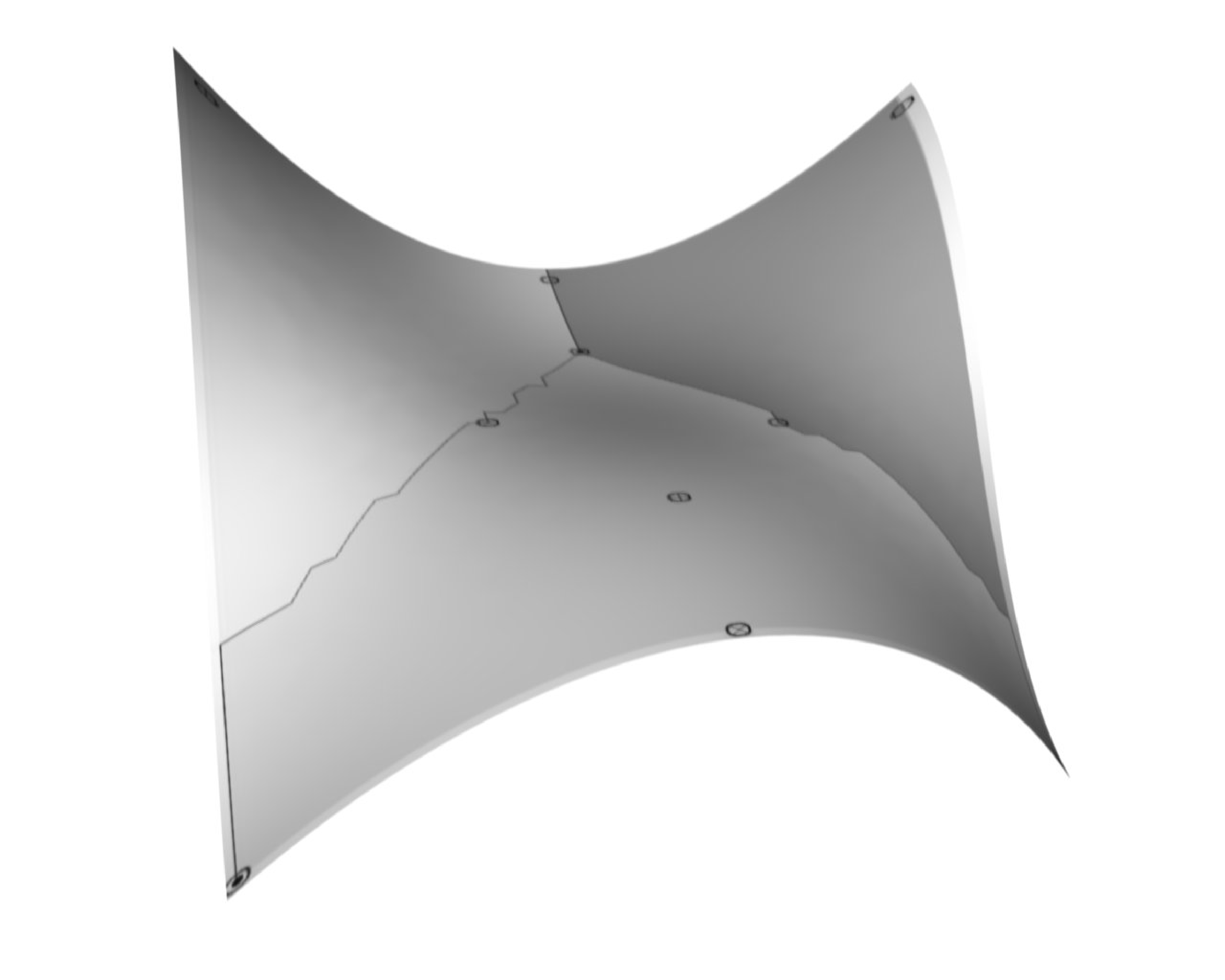}}
 \caption{Graph of the cart function with a decomposition into descending disks}
 \label{cart_discs_3d}
\end{figure}

The exact equation describing this model is
\begin{equation}\label{pospesek}
\ddot{x}=\frac{a \dot{\varphi}^2\sin(\varphi)-b\sin(2\varphi)}{c-d\cos^2\varphi}.
\end{equation}
Figure \ref{cart_discs_3d} shows the graph of this function (with suitable values of $a$, $b$, $c$ and $d$)
with the descending discs mapped onto it.

The motivation for this example comes from qualitative modeling and simulation in artificial intelligence 
\cite{Bratko03}, as a step further after the algorithm $QING$. 
The idea is, that understanding the response of the cart can provide an 
optimal strategy for controlling the system in the unstable equilibrium. 
Given a state of the system 
described by values of $\varphi$ and $\dot{\varphi}$, a $V$-path from the given state towards this point
encodes a strategy for control. Depending on some optimization criterion, an optimal strategy can be chosen. 
In the next example we give a brief explanation of a controlling strategy towards a maximum.

\subsection{Optimal path construction}
\label{teren}
\begin{figure}
\centerline{\includegraphics[width=8cm]{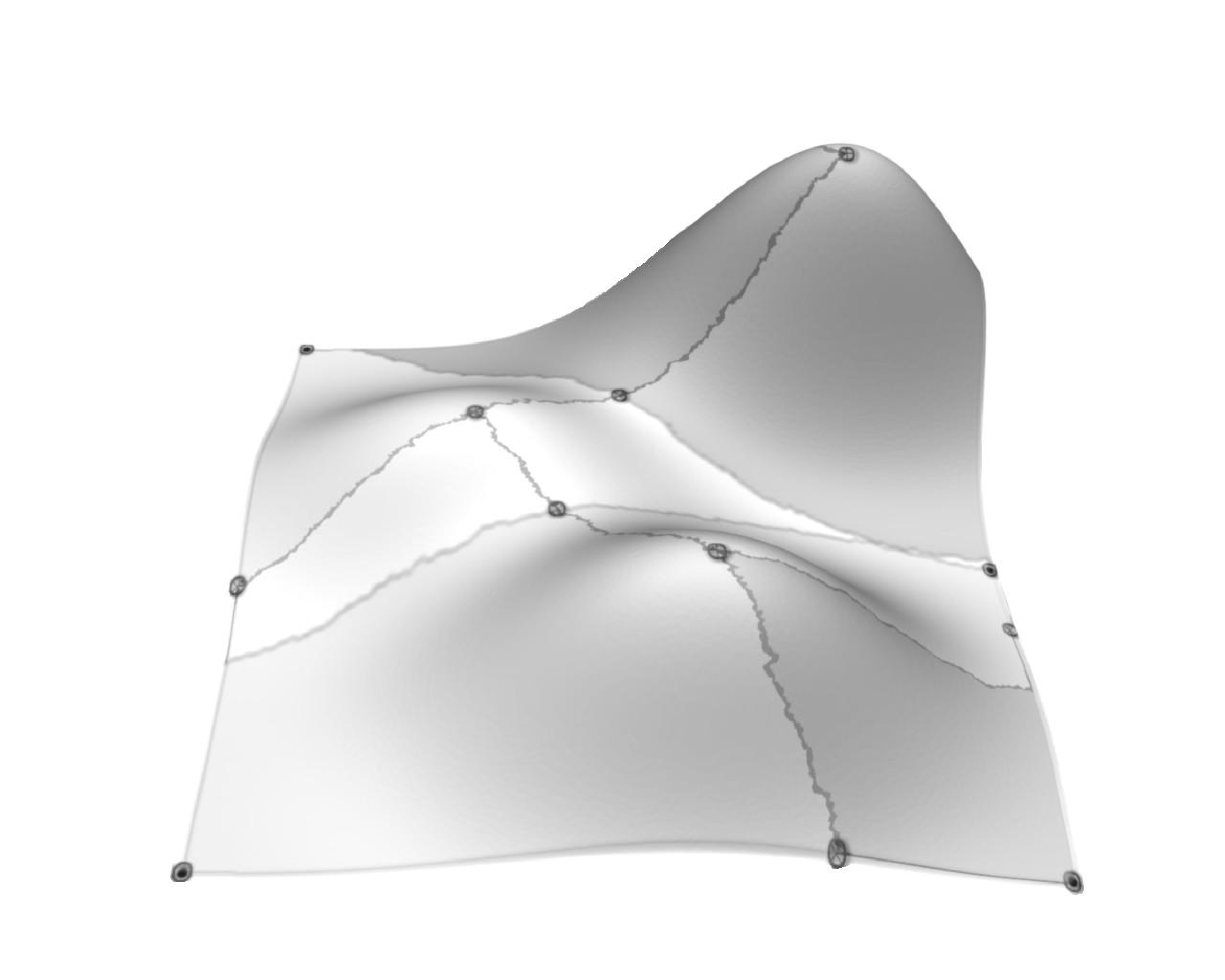}}
\caption {3D image of the terrain function on figure \ref{fig:3D_terrain}}
\label{fig:3D_terrain_decomposition}
\end{figure}

In this example, a discrete approximation of the Morse-Smale decomposition of 
the smooth function on figure \ref{fig:3D_terrain} is first constructed. The function is sampled on 
a random set consisting of $10^6$ points in the domain. 
The Delaunay triangulation on these points consists of approximately $4*10^6$ simplices. A discrete Morse 
function on this triangulation is constructed. The discrete Morse-Smale complex is shown on 
figure \ref{fig:3D_terrain_decomposition}.
The required time to build the Morse-Smale complex on a PC was a few seconds. The construction of the discrete 
Morse function was computationally the most expensive and required (on a laptop) approximately 10 minutes.  

We also present a procedure for designing an optimal path, with respect to an optimization criterion involving height 
variation, from a 
given a point $A = (x,y)$ in the domain, ending at a preferred maximal cell $\beta$ 
(for example at the highest one). A slight modification of the described procedure can be used for constructing 
a path towards any critical point or even towards any point in the domain.
Note that the procedure described can be applied in any dimension and can be therefore used to solve similar 
control problems in a more general setting. Also note, that the path constructed does not have the smallest 
possible height variation. Such a path would involve the construction of level curves which is not the object 
of this paper.

Assuming that a triangulation of the domain, a discrete Morse function on it, and a discrete Morse-Smale decomposition 
are given, we first determine the set of critical cells $S$ with the property, that $A$ lies in their descending regions. 
Since every cell is contained in at least one descending region, the set $S$ is nonempty, and because $A$ 
may belong to several descending regions, the set $S$ may contain more than one element.

Let $\alpha$ be a cell of the same dimension as $\beta$
with $a\in\alpha$. If $\simplexx{\beta}\in S$ then $\alpha$ is a cell of maximal dimension in the descending 
disk of $\beta$, so there exists a unique $V$-path from $\simplexx{\beta}$ to $\simplexx{\alpha}$, and we simply 
climb up from $\alpha$ along this path. 

If $\beta\notin S$, we construct the qualitative graph, as in paragraph \ref{qing}.
The next step is to search through paths in the graph $G$ starting from any point $\gamma\in S$ 
and ending in the required destination $\beta$ to find the path $s$ which minimizes height variation. 
Note that any other criterion can be used to choose the optimal path at this step.

The required optimal path consists of two parts. The first part connects $\alpha$ with the closest saddle  $\gamma$  
of index less than $\beta$ in $s$ which lies in the boundary of the descending region of the initial point of $s$. 
The second part of the path follows $s$ from $\gamma$ to $\beta$. 

The cell $\alpha$ is connected to $\simplexx{\gamma}$ in the following way. 
First, all $V$-paths from $\alpha$ are 
constructed. If any one of them intersects the ascending or descending disk of $\gamma$, we follow this path to the boundary, 
and then ascend or descend to $\gamma$. If none of them do, we make a step along the unique $V$-path from the initial point of $s$ 
to $\alpha$, and repeat the procedure. After finitely many steps we either reach the descending or ascending 
region of  the initial point of $s$, or $\alpha$, from where a $V$-path to this point exists. 

\subsection{A macroeconomic model}

In our last application, macroeconomic financial indicators of European countries 
are analyzed. The data was retrieved from the publicly available database Eurostat \cite{eurostat}. This example is part of 
an ongoing project involving  
a qualitative analysis of the effect of key macroeconomic indicators on long term economic growth of a country.  
The motivation for this project is the hypothesis 
formulated in \cite{DMVB} that long term economic growth can be achieved by balancing key macroeconomic factors.  
In this example a five-year average of real GDP (gross domestic product) growth rate $(Ggr)$ is used as 
a measure of long term economic growth, and the effect of the following macroeconomic indicators as 
independent variables is analyzed: inflation rate ($I$), balance of current account ($BCA$), and public debt ($PD$). 
An additional independent variable, GDP per capita ($Gpp$), was used as a control variable to distinguish between 
different levels of economic development in the countries. 

Each point in our model thus represents a point in $\mathbb{R}^4$ with coordinates the 
values of $I$, $BCA$ and $PD$ for one of the European countries in a given year from 1998 to 2004, 
and the average value of $Gpp$ over this and the next four 
years. The average value of $Ggr$ over this and the next four years is the dependent variable. 
Altogether 143 data points representing the  
countries included in different years were available. 

\begin{figure}[ht!]
  \centerline{\includegraphics[width=14cm]{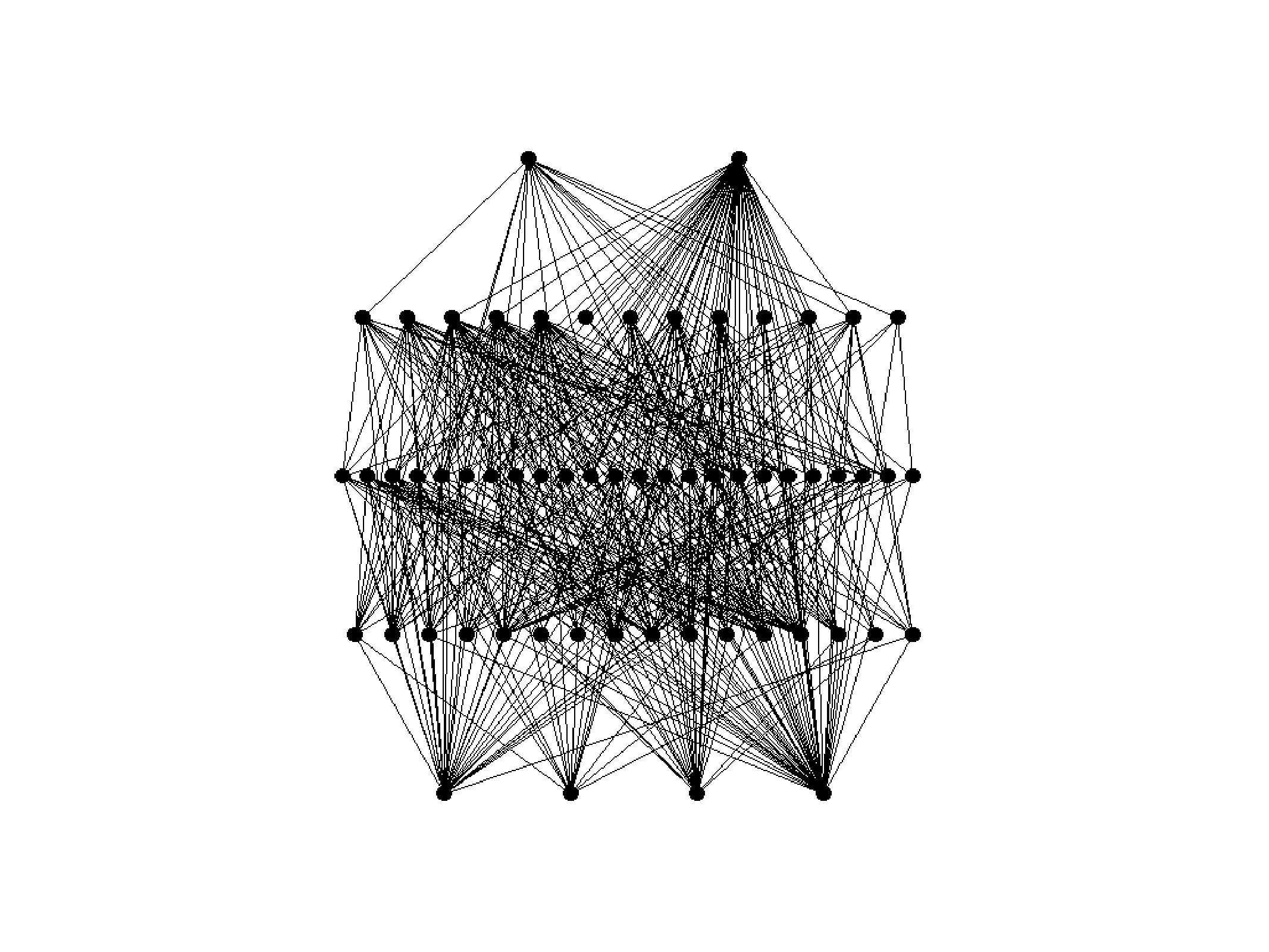}}
  \label{slika_grafa}
  \caption{The qualitative graph connecting the critical cells in the 4-dimensional macroeconomic data set.}
\end{figure}

A Delaunay triangulation on the points in $\mathbb{R}^4$ was constructed, producing altogether 17451 simplices of dimension up to 
4. The discrete vector field was reconstructed from the data points using the algorithm of \cite{KKM05} 
producing 2 critical cells in dimension 4 (maxima), 13 in dimension 3, 24 in dimension 2, 16 in dimension 1 
and 4 critical cells in dimension 0. On a laptop, the descending and ascending disks are reconstructed immediately, 
once the triangulation and the discrete vector field are given. The reconstruction of the triangulation is also almost 
immediate, and the discrete gradient vector field is built in less than a minute.  
The graph, connecting the critical cells with respect to 
incidence in the boundaries of the corresponding descending regions is given in figure \ref{slika_grafa}.

For each country, data for the current year determines a point in the domain, and the model can be used to 
predict its average economic growth in the next five years. 
The $V$-paths encode the changes in the parameter values 
which ensure constant growth of the GDP growth rate. Finally, the Morse-Smale decomposition can be used 
to construct a path, similar to the one in example \ref{teren}, which ends in a preferred position 
involving minimal variation in the GDP growth.

The two maxima correspond to Ireland in 1999 ($m_1$) and Latvia in 2003 ($m_2$). 
According expert knowledge from economics, $m_1$ represents the preferred maximum
(since the high growth in Latvia is 
partly due to its post transition stage from a state planned economy, and involves high risks in 
the form of very high inflation, public 
debt and unemployment). So, for a country which is currently in the descending disc of $m_2$ 
(for example Slovenia), 
a path leading from its current position to a common 
$3$-saddle, and 
from there along the gradient vector field to $m_1$ encodes the recommended strategy.

%\bibliography{gregor}{}
%\bibliographystyle{elsart-num}
%\bibliographystyle{splncs}
% \bibliographystyle{plain}

\end{document}